\documentclass[11pt]{amsart}
\usepackage{amsmath,amssymb,amscd,graphicx,color}

\newtheorem{thm}{Theorem}[section]

\newtheorem{prop}[thm]{Proposition}

\theoremstyle{definition}
\newtheorem{defn}[thm]{Definition}
\newtheorem{rem}[thm]{Remark}
\newtheorem{exmp}[thm]{Example}

\newcommand{\Z}{\ensuremath{\mathbb{Z}}}
\newcommand{\D}{\ensuremath{\mathbb{D}}}
\newcommand{\R}{\ensuremath{\mathbb{R}}}
\newcommand{\F}{\ensuremath{\mathbb{F}}}

\newcommand{\sph}{\ensuremath{\mathbb{S}}}
\newcommand{\smallcaps}[1]{\ensuremath{\textsc{#1}}}
\newcommand{\tor}{\smallcaps{Tor}}
\newcommand{\art}{\smallcaps{Art}}
\newcommand{\script}[1]{\ensuremath{\mathcal{#1}}}

\newcommand{\rel}{\script{R}}

\newcommand{\onto}{\ensuremath{\twoheadrightarrow}}

\DeclareMathOperator{\link}{Link}
\newcommand{\bt}{\begin{tabular}}
\newcommand{\et}{\end{tabular}}
\newcommand{\ba}{\begin{array}}
\newcommand{\ea}{\end{array}}

\begin{document}
%%%%%%%%%%%%%%%%

\title[Multi-vertex complexes]
  {Combinatorial descriptions of \\ multi-vertex $2$-complexes}
\author[Jon McCammond]{Jon McCammond$\ \!{ }^1$}
      \address{Dept. of Math.\\
               UC Santa Barbara\\
               Santa Barbara, CA 93106}
      \email{jon.mccammond@math.ucsb.edu}
\date{\today}

\begin{abstract}
  Group presentations are implicit descriptions of $2$-dimensional
  cell complexes with only one vertex.  While such complexes are
  usually sufficient for topological investigations of groups,
  multi-vertex complexes are often preferable when the focus shifts to
  geometric considerations.  In this article, I show how to quickly
  describe the most important multi-vertex $2$-complexes using a
  slight variation of the traditional group presentation.  As an
  illustration I describe multi-vertex $2$-complexes for torus knot
  groups and one-relator Artin groups from which their elementary
  properties are easily derived.  The latter are used to give
  an easy geometric proof of a classic result of Appel and Schupp.
\end{abstract}
\maketitle

\footnotetext[1]{Partially supported by grants from the National
  Science Foundation}

Some cell complexes are easy to describe: a graph with one vertex
corresponds to a set $S$ indexing its edges and a one-vertex
combinatorial $2$-complex can be constructed from an algebraic
presentation $\langle S \mid \rel \rangle$.  When one tries to
describe $2$-complexes with multiple vertices, however, several issues
arise.  First, there is no standard way to quickly describe a
complicated $1$-complex. And second, even supposing such a
$1$-skeleton as given with edges oriented and labeled by a set $S$,
not all words over the alphabet $S \cup S^{-1}$ can be used to
describe closed paths, making it easy to list collections of words
that are incompatible with the given graph.  In this article we
describe a simple procedure that avoids both of these difficulties and
requires only mild restrictions.  It constructs a multi-vertex
link-connected combinatorial $2$-complex from any multiset of words,
and every such complex can be constructed in this way.  Such a process
is sufficient for most purposes since the only $2$-complexes excluded
are those that are homotopy equivalent to a non-trivial wedge product,
i.e. those whose fundamental groups can be freely decomposed.  After
describing this procedure and establishing its main properties, sample
applications are given that illustrate how multi-vertex complexes can
make geometric properties of groups more transparent, including a
short geometric proof of a classic result of Appel and Schupp
\cite{ApSch83}.  

%On a final note, it is with great pleasure that I
%dedicate this article to Paul Schupp for all of his wonderful
%mathematics and for his inspiration over the years.

%%%%%%%%%%%%%%%%%%%%%%%%%%%%%%%%%%%%%%
\section{Standard $2$-complexes}
%%%%%%%%%%%%%%%%%%%%%%%%%%%%%%%%%%%%%%

We begin by reviewing the standard method for creating a one-vertex
combinatorial $2$-complex from an algebraic presentation $\langle S
\mid \rel \rangle$.  Recall that CW complexes are inductively
constructed by attaching $n$-discs along their boundary cycles to an
already constructed $(n-1)$-skeleton and that $1$-dimensional CW
complexes are undirected graphs.  Recall also that a map $Y\rightarrow
X$ between CW complexes is a \emph{combinatorial map} if its
restriction to each open cell of $Y$ is a homeomorphism onto an open
cell of $X$ and that a CW complex $X$ is \emph{combinatorial} provided
that the attaching map of each cell of $X$ is combinatorial for a
suitable subdivision of its domain.  In this article all maps and cell
complexes are combinatorial unless otherwise specified.  For a
$2$-complex, this means that it can be viewed as the result of
attahing polygons to a graph using combinatorial maps.

\begin{defn}[Polygons and $2$-complexes]\label{def:2-complexes}
  A \emph{polygon} is a (closed) $2$-disc $\D^2$ whose boundary cycle
  has been given the structure of a graph.  When its boundary cycle
  has combinatorial length $n$ it is called an \emph{$n$-gon}.  A
  $2$-complex $X$ is constructed from a graph $\Gamma$ and a
  collection $\script{P}$ of disjoint polygons by specifying for each
  $n$-gon in $\script{P}$ a closed combinatorial path of length~$n$ in
  $\Gamma$ along which its boundary cycle should be attached.  Let $P$
  denote the disjoint union of polygons in $\script{P}$ and note that
  $P$ itself is a cell complex.  Also note that so long as $X$ has no
  isolated vertices and every edge of $\Gamma$ occurs in the image of
  at least one attaching map, the complex $X$ is a quotient of the
  complex $P$ and the quotient map $P \onto X$ is a combinatorial map.
  When $X$ satisfies these minor restrictons, we say that $X$ is a
  \emph{polygon quotient} with quotient map $P \onto X$.
\end{defn}

A $2$-complex with only one vertex is called a \emph{standard
  $2$-complex} and the traditional way to quickly and efficiently
describe it is via a presentation.

\begin{defn}[Presentations]\label{def:presentations}
  A \emph{presentation} $\langle S \mid \rel \rangle$ consists of a
  set $S$ and a multiset of words $\rel$ over the alphabet $S \cup
  S^{-1}$.  (We say \emph{multiset} rather than set because
  repetitions are allowed.)  The elements of $S$ are \emph{generators}
  and the elements of $\rel$ are \emph{relators}.
\end{defn}

Presentations and standard $2$-complexes are essentially
interchangeable.

\begin{thm}[Standard $2$-complexes]\label{thm:std-2-complexes}
  Every presentation $\langle S \mid \rel \rangle$ implicitly
  descrpibes a standard $2$-complex and every standard $2$-complex can
  be constructed from a presentation.
\end{thm}

\begin{proof}
  To construct a standard $2$-complex from a presentation $\langle S
  \mid \rel \rangle$ first use the set $S$ to build a one-vertex
  directed graph $\Gamma$ whose edges are indexed by $S$ and note that
  combinatorial paths in $\Gamma$ are in natural bijection with words
  over the alphabet $S \cup S^{-1}$.  Next, for each word of
  length~$n$ in $\rel$, we attach an $n$-gon to $\Gamma$, identifying
  its (based and oriented) boundary cycle with the combinatorial path
  of length~$n$ in $\Gamma$ that corresponds to this word.

  In the other direction, given a standard $2$-complex $X$ with
  $1$-skeleton $\Gamma$, one chooses orientations for the edges of
  $\Gamma$ and indexes them by a set $S$.  Then, for each polygon
  attached to $\Gamma$, choose a base vertex in its boundary and an
  orientation of its boundary cycle.  The attaching map of this
  $2$-cell can then be encoded in the word corresponding to the closed
  combinatorial path in $\Gamma$ described by the image of this based
  and oriented boundary cycle.  If $\rel$ denotes the collection of
  these words, it should be clear that the presentation $\langle S
  \mid \rel \rangle$ can be used to reconstruct the standard
  $2$-complex $X$.
\end{proof}

When describing concrete examples we use several simplifying
conventions.  Uppercase roman letters are used to denote the inverse
of their lowercase equivalents in order to make words easier to parse
and absorb.  Thus, we write $abAB$ instead of $aba^{-1}b^{-1}$.  We
also allow relators to be given implicitly via relations.  A
\emph{relation} is an equation of the form $r = s$ where $r$ and $s$
are words over the alphabet $S \cup S^{-1}$ and the implicit relator
is the word $rs^{-1}$.  For example, the relation $ab=ba$ refers
implicitly to the relator $abAB$.

%%%%%%%%%%%%%%%%%%%%%%%%%%%%%%%%%%%%%%%%%
\section{Multi-vertex $2$-complexes}
%%%%%%%%%%%%%%%%%%%%%%%%%%%%%%%%%%%%%%%%%

In this section we introduce an alternative construction.

\begin{defn}[Constructed by edge identifications]\label{def:edge-ident}
  Let $X$ be a $2$-complex that is a polygon quotient and let $P \onto
  X$ be the corresponding quotient map
  (Definition~\ref{def:2-complexes}).  A third cell complex $Y$,
  between $P$ and $X$, can be defined as follows.  Identify pairs of
  $1$-cells in $P$ iff they are sent to the same $1$-cell in $X$, and
  identify them in the same fashion.  For $Y$ to be a cell complex
  certain vertex identifications must also be made, but make only
  those identifications that are forced by the edge identifications.
  The polygon quotient map $P \onto X$ thus factors into two
  combinatorial maps $P \onto Y \onto X$ and we say that $Y$ is
  \emph{constructed from $X$ by edge identifications}.  Finally, note
  that since $P \onto Y$ is a factor of $P \onto X$, the only vertices
  in $P$ that can be identified in $Y$ are those with the same image
  in $X$.
\end{defn}  

In order to clarify under what conditions the map $Y \onto X$ is a
homeomorphism, we recall the notion of a vertex link.

\begin{defn}[Vertex links]
  Let $X$ be a polygon quotient with quotient map $f:P \onto X$.  For
  each vertex $u$ in $X$ there is a $1$-complex $\link(u,X)$ called
  the \emph{link of $u$ in $X$}.  Intuitively, it is the boundary of
  an $\epsilon$-neighborhood around $u$ in $X$, but in the absence of
  a metric, one can also define it as follows.  Start with a distinct
  closed edge for each vertex $v$ in $P$ and associate its two
  endpoints with the two ends of edges attached to $v$.  Next,
  restrict attention to those closed edges associated with vertices
  $v$ with $f(v)=u$.  Finally, identify the endpoints of the closed
  edges iff the corresponding ends of edges in $P$ are identified
  under the quotient map $f$.  The result is the graph $\link(u,X)$.
  We say that $X$ is a \emph{link-connected $2$-complex} when for
  every vertex $u$ in $X$, $\link(u,X)$ is a connected graph.
\end{defn}

\begin{prop}[Identifying vertices]\label{prop:vertex-ident}
  Let $X$ be a polygon quotient with quotient map $f:P \onto X$ and
  let $Y$ be the complex constructed from $X$ by edge identifications.
  If $v$ and $v'$ are vertices in $P$ with $f(v) = f(v') = u$ in $X$,
  then $v$ and $v'$ are identified in $Y$ iff the edges of
  $\link(u,X)$ corresponding to $v$ and $v'$ belong to the same
  connected component.  As a consequence $Y$ is always link-connected
  and the map $Y\onto X$ is a homeomorphism iff $X$ itself is
  link-connected.
\end{prop}

\begin{proof}
  Both directions of the first assertion are straightforward.  If the
  corresponding edges belong to the same connected component then
  there is a finite length path connecting them in the link and this
  path encodes a finite sequence of individual edge identifications
  that force $v$ and $v'$ to be identified in $Y$.  Conversely,
  identifying vertices iff the corresponding edges belong to the same
  connected component of the link produces an intermediate cell
  complex in which all the edge identifications can be performed with
  no further vertex identifications.  Thus, no additional vertex
  identifications are forced.

  Next note as a consequence of the first assertion that the vertex
  links of $Y$ are connected components of vertex links of $X$.  Thus
  $Y$ is link-connected.  Moreover, when $X$ is link-connected, the
  map $Y\onto X$ is bijective on vertices and an isomorphism on vertex
  links and it quickly follows that it is a homeomorphism.
  Conversely, if $X$ has a single disconnected vertex link, then
  distinct vertices of $Y$ are identified in $X$ and the map $Y \onto
  X$ not a homeomorphism.
\end{proof}

Now that these properties have been established, we turn our attention
to constructing a link-connected $2$-complex from a multiset of words.

\begin{defn}[Combinatorial descriptions]\label{def:comb-descr}
  Let $S$ be a set and let $\rel$ be a nonempty multiset of words over
  the alphabet $S \cup S^{-1}$.  The list $[\rel]$ is called a
  \emph{combinatorial description} and square brackets are used to
  instead of angle brackets to highlight that this is not a
  traditional presentation. The elements of $\rel$ are still called
  relators and the same simplifying conventions remain in effect.
\end{defn}

Our main result is that combinatorial descriptions and link-connected
$2$-complexes are essentially interchangeable.

\begin{thm}[Link-connected $2$-complexes]\label{thm:link-conn-2-complex}
  Every combinatorial description $[\rel]$ implicitly describes a
  link-connected $2$-complex $Y$ and every link-connected $2$-complex
  $Y$ can be constructed from a combinatorial description.
\end{thm}

\begin{proof}
  Let $[\rel]$ be a combinatorial description, let $S$ be the set of
  letters that occur in the relators in $\rel$, and let $X$ be the
  standard $2$-complex described by the presentation $\langle S \mid
  \rel \rangle$.  Because of the restriction on $S$, $X$ is a polygon
  quotient and we can define $Y$ as the complex constructed from $X$
  by edge identifications.  By Proposition~\ref{prop:vertex-ident} $Y$
  is link-connected.  Alternatively, and more directly, we can proceed
  as follows.  First, let $P$ be a disjoint union of polygons indexed
  by the words in $\rel$ so that words of length~$n$ in $\rel$
  correspond to $n$-gons.  Next, orient and label the edges in the
  boundary cycle of each polygon according to its corresponding word.
  (Using the standard $2$-complex $X$ this can be done by pulling back
  the labels and orientations of the edges in the one-vertex graph
  $\Gamma$ derived from $S$ through the attaching maps of the
  $2$-cells of $X$.)  Finally, rather than using the labeled oriented
  edges of $P$ to identify how these boundary cycles should be
  attached to $\Gamma$, we use this information instead to identify
  which of these edges should be identified with each other.  In
  particular, $Y$ is the quotient of $P$ which identifies edges
  according to label and orientation, and which identifies vertices
  iff the identification is necessary so that the quotient remains a
  cell complex.
  
  In the other direction, given a link-connected $2$-complex $Y$, with
  $1$-skeleton $\Gamma$, one chooses orientations for the edges of
  $\Gamma$ and indexes them by a set $S$.  Then, for each polygon
  attached to $\Gamma$, choose a base vertex in its boundary and an
  orientation of its boundary cycle.  The attaching map of this
  $2$-cell can then be encoded in the word corresponding to the closed
  combinatorial path in $\Gamma$ described by the image of this based
  and oriented boundary cycle.  If $\rel$ denotes the collection of
  these words, it should be clear that the combinatorial description
  $[\rel]$ can be used to reconstruct a link-connected $2$-complex
  that is equal to $Y$ by Proposition~\ref{prop:vertex-ident}.
\end{proof}

When a combinatorial description $[\rel]$ and a link-connected
$2$-complex $Y$ are related in this way we say $Y$ is the
\emph{complex constructed from $[\rel]$} and $[\rel]$ is a
\emph{combinatorial description of $Y$}.  Note that we use the word
``combinatorial'' rather than ``algebraic'' since the letters in $S$
correspond to edges with possibly distinct endpoints.  In particular
they need not be closed loops and thus do not have a natural algebraic
intepretation in the fundamental group of $Y$.  The distinction
between combinatorial descriptions and presentations is highlighted by
the following example.

\begin{exmp}[Descriptions vs. presentations]\label{ex:desc-vs-pres}
  The $2$-complex constructed from the combinatorial description
  $[abcABC]$ (or equivalently $[abc=cba]$) is a torus with two
  vertices and thus its fundamental group is $\Z^2$.  (More generally,
  any combinatorial description in which every letter occurs exactly
  twice--in either orientation--corresponds to a closed surface.)  The
  presentation $\langle a,b,c \mid abc=cba \rangle$, on the other
  hand, corresponds to a quotient of this torus with its two vertices
  identified.  Since it is homotopy equivalent to a wedge product of a
  torus and a circle, its fundamental group is $\Z^2 \ast \Z$.
\end{exmp}

%%%%%%%%%%%%%%%%%%%%%%%%%%%%
\section{Wedge products}
%%%%%%%%%%%%%%%%%%%%%%%%%%%%

Standard $2$-complexes are considered sufficiently flexible for most
purposes since every connected $2$-complex $X$ is homotopy equivalent
to a standard $2$-complex; one simply selects a spanning tree in the
$1$-skeleton of $X$ and collapses it to a point.  In this section we
show that link-connected $2$-complexes are nearly as flexible by
establishing the following result.

\begin{thm}[Splitting $2$-complexes]\label{thm:splitting}
  Every group is the fundamental group of a wedge product of circles
  and link-connected $2$-complexes.
\end{thm}

\begin{figure}
\includegraphics[width=4in]{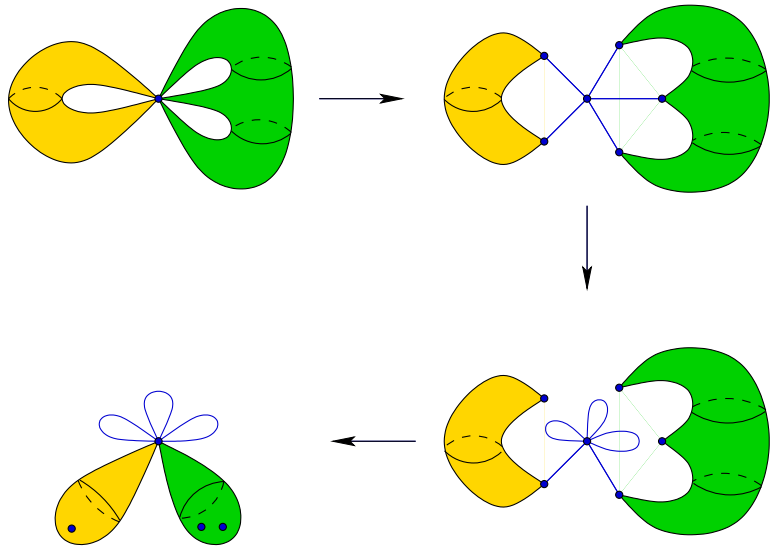}
\caption{An illustration of the homotopy equivalences used to convert
 an arbitrary $2$-complex into a wedge product of circles and
 link-connected $2$-complexes.\label{fig:splitting}}
\end{figure}

\begin{proof}
  Let $G$ be a group, let $X$ be a standard $2$-complex with $G$ as
  its fundamental group, and let $L=\link(\ast,X)$ where $\ast$ is the
  unique vertex of $X$.  The proof proceeds by repeatedly modifying
  $X$ using a series of homotopy equivalences.  An illustration of the
  process is shown in Figure~\ref{fig:splitting}.  When the link $L$
  is connected, there is nothing to prove, so suppose not and let $I$
  and $J$ be sets that index the connected components of the link $L$
  and the connected components of $X\setminus \{\ast\}$,
  respectively. Note that since the link $L$ can be viewed as the
  boundary of an $\epsilon$-neighborhood of $\ast$ in $X$, there is a
  well-defined function $f:I\onto J$.

  We construct a new $2$-complex $Y$ by pulling the connected
  components of $L$ in different directions.  More specifically, start
  with a tree $T$ that has $0$-cells indexed by $I \sqcup \{\ast\}$
  and an edge $e_i$ from $v_\ast$ to $v_i$ for each $i \in I$.  The
  rest of $Y$ is built by adding a $1$-cell or $2$-cell to $T$ for
  each $1$-cell and $2$-cell in $X$ in such a way that the complex
  obtained by contracting $T$ to a point is equal to $X$.  Concretely,
  for each $1$-cell of $X$ we add a $1$-cell to $T$ with each end
  attached to the vertex $v_i$ in $T$ where $i\in I$ indexes the
  component of $L$ through which this end approaches $\ast$ in $X$.
  This completes the $1$-skeleton.

  For each $2$-cell of $X$ we attach a $2$-cell to $Y^{(1)}$ along the
  corresponding sequence of edges.  Because of the way the edges of
  $X$ were attached to $T$, the old closed combinatorial paths in the
  $1$-skeleton of $X$ correspond to closed combinatorial paths in the
  $1$-skeleton of $Y$.  More specifically, because paths of length~$2$
  in the boundary cycles of $2$-cells create edges in $L$, the ends of
  these adjacent edges belong to the same component $i$, their lifts
  are attached to the same vertex $v_i$, and thus the new edges can be
  concatenated as before.  Since collapsing the contractible
  subcomplex $T$ to a point converts $Y$ into $X$, the two are
  homotopy equivalent.
  
  The remaining steps are similarly straightforward.  Since
  $Y\setminus T$ is homeomorphic to $X\setminus \{\ast\}$ under the
  quotient map, its connected components remain indexed by $J$.  For
  each $j\in J$ select an edge $e_i$ with $f(i)=j$ and then reattach
  all unselected edges in $T$ so that both of their endpoints are at
  $v_\ast$.  See the lower righthand corner of
  Figure~\ref{fig:splitting}.  The result is homotopy equivalent to
  $Y$ since there is a path from the other endpoint to $v_\ast$ that
  travels through a component of $Y\setminus T$ and then back to
  $v_\ast$ along a selected edge, making the original and altered
  attaching maps homotopic.

  The last step is to contract the tree formed by the selected edges
  to a point and to note that the result is a wedge product of circles
  and complexes indexed by $J$. Every vertex link in a complex indexed
  by $J$ is connected since, by construction, it can be identified
  with a connected component of the original link $L$.
\end{proof}

A corollary of Theorem~\ref{thm:splitting} is that every group that
does split as a non-trivial free product is the fundamental group of a
link-connected $2$-complex.  We conclude this short section with a
concrete illustration of the proof.

\begin{exmp}[Splitting $2$-complexes]\label{ex:splitting}
  Let $X$ be the quotient of two disjoint $2$-spheres that identifies
  two distinct points in the first $2$-sphere and three distinct
  points in the second $2$-sphere to a single point.  The quotient $X$
  can be given a cell structure so that it is a standard $2$-complex,
  but the exact cell structure is irrelevant.  The link of the unique
  vertex $\ast$ in $X$ has $5$ connected components and $X\setminus
  \{\ast\}$ has $2$.  In other words $|I|=5$ and $|J|=2$.
  Figure~\ref{fig:splitting} illustrates the sequence of steps used to
  show that $X$ is homotopy equivalent to $\sph^2 \vee \sph^2 \vee
  \sph^1 \vee \sph^1 \vee \sph^1$.
\end{exmp}

%%%%%%%%%%%%%%%%%%%%%%%%
\section{Torus knots}
%%%%%%%%%%%%%%%%%%%%%%%%

Having established properties of link-connected $2$-complexes and
their combinatorial descriptions, we now turn our attention to
examples that illustrate the benefits of using multi-vertex
$2$-complexes.  We begin with torus knots and torus knot groups.

\begin{defn}[Torus knots]
  The $3$-sphere has a standard genus one Heegaard splitting into two
  solid tori with a common torus boundary and any simple closed curve
  that embeds in this common torus is called a \emph{torus knot}.  The
  essential curves on this torus that bound discs in one solid torus
  or the other provide a canonical basis for the first homology of the
  torus and torus knots can be classified by the element of first
  homology they represent.  In particular, for every relatively prime
  pair of integers $p$ and $q$ there is a knot $K$ called a
  \emph{$(p,q)$-torus knot} corresponding to a $(p,q)$-curve on this
  separating torus.  The fundamental group of the complement of $K$ is
  the corresponding \emph{torus knot group} $\tor(p,q)$ and a
  presentation of this group is $\langle a,b \mid a^p = b^q \rangle$.
  Although torus knots are only defined when $p$ and $q$ be relatively
  prime, the presentation makes sense for arbitrary pairs of integers
  $m$ and $n$ and we extend the definition of $\tor(m,n)$ accoringly.
\end{defn}

\begin{figure}
  \bt{cc}
  \bt{c}\includegraphics{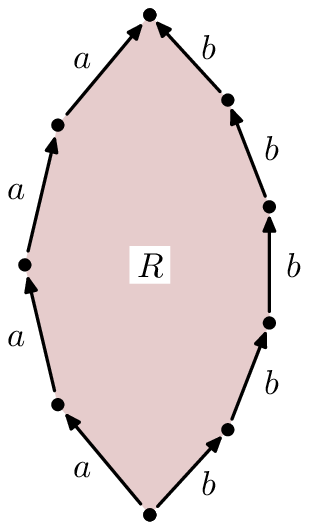}\et & 
  \bt{c}\includegraphics{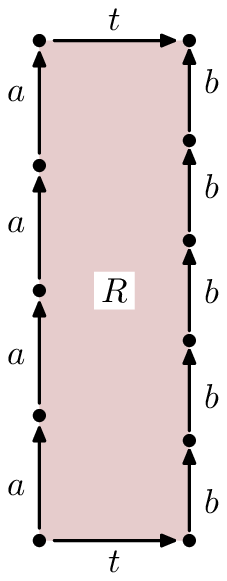}\et \\
  \bt{c}\includegraphics{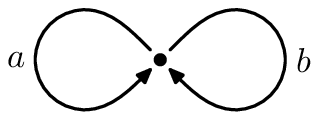}\et & 
  \bt{c}\includegraphics{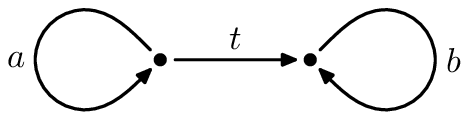}\et \\
  \et  
  \caption{The $2$-cell on the left attached to the graph on the left
    and the $2$-cell on the right attached to the graph on the right
    are both $2$-complexes with fundamental group $\tor(4,5)$.
    Contracting the edge $t$ shows that the two spaces are homotopy
    equivalent.\label{fig:torus-cells}}
\end{figure}

Let $X = X_{m,n}$ be the standard $2$-complex of the standard
presentation of $\tor(m,n)$, i.e. $\langle a,b\mid a^m = b^n\rangle$.
Although the presentation of a torus knot group is extremely simple,
the global structure of the universal cover of $X$ is not immediately
obvious.  The situation is much clearer if we consider the two vertex
$2$-complex $Y = Y_{m,n}$ corresponding to the combinatorial
description $[a^mt = tb^n]$.  The $2$-cells and the underlying graphs
of the complexes $X_{4,5}$ and $Y_{4,5}$ are depicted in
Figure~\ref{fig:torus-cells}.  That $X$ and $Y$ are two $2$-complexes
with the same fundamental group is clear since the edge labeled $t$ in
$Y$ is always embedded and contracting it to a point yields $X$.  

\begin{figure}
\includegraphics{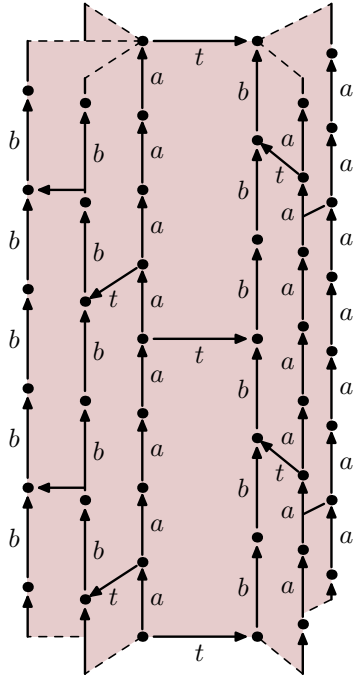}
\caption{A portion of the universal cover of $Y_{4,3}$.\label{fig:torus-cover}}
\end{figure}

As mentioned above, the main benefit of using $Y$ instead of $X$ as a
space with fundmamental group $\tor(m,n)$ is that the universal cover
of $Y$ is much easier to visualize in its entirety.  First add a
metric to the polygon used to construct $Y$.  As is hinted in
Figure~\ref{fig:torus-cells}, we turn it into a metric rectangle with
right angles at the four endpoints of the two edges labeled $t$ and
with no other sharp corners.  To make the edge lengths match up, we
make the $a$ edges length~$n$ and the $b$ edges length~$m$.  In the
universal cover, these rectangles glue together along the $t$-edges to
produce vertical strips that are in turn glued together in a tree-like
fashion.  In fact, the universal cover $\widetilde{Y}$ can be
described as a metric direct product of a tree $T$ and a copy of the
real line $\R$.  See Figure~\ref{fig:torus-cover}.  If we let
$K=K_{m,n}$ denote the complete bipartite graph with $m$ vertices of
one type and $n$ vertices of the other type, then the tree $T$ is the
universal cover of $K$.  In particular, $T$ is \emph{biregular} in
that every vertex has valence $m$ or $n$ and every edge connects a
vertex of valence $m$ to a vertex of valence $n$.  As usual with
covering spaces, the group $\tor(m,n)$ acts freely and cocompactly by
isometries on the metric space $\widetilde Y = T \times \R$, which is
contractible and non-positively curved.  Using the action of
$\tor(m,n)$ on this space it is straightforward to establish the
following elementary properties of torus knot groups.

\begin{thm}[Torus knot groups]\label{thm:torus}
  If $G = \tor(m,n)$ is a torus knot group for positive integers $m$
  and $n$ then: (1) the center of $G$ is infinite cyclic generated by
  the element $a^m = b^n$; (2) $G$ is virtually a direct product of a
  free group and an infinite cyclic group; (3) every nontrivial
  reduced word equivalent to the identity in $G$ contains a subword
  equal to $a^m$ or $b^n$ or their inverses; and finally, (4) every
  word equivalent to the identity in $G$ can be reduced to the
  identity by iteratively replacing $a^m$ with $b^n$, replacing $b^n$
  with $a^m$, and performing free reductions.
\end{thm}

\begin{proof}[Proof Sketch]
  Since the listed properties are relatively elementary and they
  follow fairly quickly once the geometry of $\widetilde Y$ and action
  of $G$ is understood, we merely sketch the proofs.  First note that
  because contracting the edge labeled $t$ in $Y$ yields $X$,
  contracting the disjoint edges labeled $t$ in $\widetilde Y$ yields
  $\widetilde X$.  Thus, if we treat the edges labeled $t$ in
  $\widetilde Y$ as though they were contracted (without actually
  contracting them) we can work with the geometrically pleasing
  $1$-skeleton of $\widetilde Y$ to establish results about the
  $1$-skeleton of $\widetilde X$, i.e. the Cayley graph of $G$ with
  respect to the generating set $\{a,b\}$.  For example, in $X$, the
  edges labeled $a$ and $b$ label loops which represent elements of
  the fundamental group, and as such they act on $\widetilde X$ by
  deck transformations once we have chosen a vertex in $\widetilde X$
  as our base vertex.  Thus they also represent actions on $\widetilde
  Y$ by deck transformations once we have chosen an edge labeled $t$
  as our base edge.

  Fix such an edge $t$ and consider the column labeled with $a$'s at
  one endpoint, the column labeled with $b$'s at the other and the
  vertical strip between them.  The deck transformation corresponding
  to the generator $a$ shifts this $a$ column vertically and spins the
  rest of $\widetilde Y$ around this column.  After $m$ such motions
  the entire complex $\widetilde Y$ merely experiences a vertical
  shift.  Similarly, the deck transformation corresponding to the
  generator $b$ shifts this $b$ column vertically and spins the rest
  of $\widetilde Y$ around this column.  After $n$ such motions the
  entire complex $\widetilde Y$ merely experiences a vertical shift.
  From these motions one can show that the only words that commute
  with $a$ correspond to paths that start at the basic $t$ edge and
  end at another $t$ edge with an endpoint on this $a$ column.
  Similarly, any word that commutes with $b$ must correspond to a path
  that starts at the basic $t$ edge and ends at another $t$ edge with
  an endpoint on this $b$ column.  Thus, the only words that might be
  central are those that correspond to a path that starts at the basic
  $t$ edge and ends at another $t$ edge in the same vertical strip.
  As all these words represent rigid vertical shifts and are powers of
  the basic vertical shift represented by $a^m=b^n$, the infinite
  cyclic subgroup generated by this element is precisely the center of
  $G$.

  To see (2) we note that there is a finite-sheeted cover of $Y$
  obtained by identifying edges labeled $t$ in $\widetilde Y$ when
  they belong to the same vertical strip and also identifying two
  vertical strips if they have their $t$ edges at the same set of
  heights.  We then make the minimal additional identifications
  necessary for the result to be a covering of $Y$.  Geometrically the
  result is a direct product $K_{m,n} \times \sph^1$ with fundamental
  group $\F \times \Z$ where $\F$ is the free group $\pi_1(K_{m,n})$
  and since the cover is finite-sheeted, the subgroup this represents
  is finite index.

  Next, recall that a \emph{syllable} of a word is a maximum subword
  that merely repeats the same letter.  For example, the word
  $a^5b^2C^4$ has $3$ syllables: $a^5$, $b^2$, and $C^4$
  (i.e. $c^{-4}$).  For (3) we convert a reduced word equivalent to
  the identity into a closed immeresed path in the $1$-skeleton of
  $\widetilde X$ starting at its base vertex and then finally to a
  closed immersed path in the $1$-skeleton of $\widetilde Y$ by
  traversing $t$ edges when necessary in order to continue (which
  occurs precisely at the breaks between syllables).  Given such a
  path, we can look at its projection into the tree $T$.  The
  projection cannot be trivial since there are no closed immersed
  paths that remain in a single $a$ column or a single $b$ column.
  Also, the projected curve cannot remain immersed since $T$ is a
  tree.  Thus there is a point in the projected curve where it crosses
  an edge of $T$ and then immediately backtracks across the same edge.
  If we consider the portion of the path in $\widetilde Y$ that
  produces this behavior, we see a path that crosses a $t$ edge,
  travels up or down an $a$ column (or $b$ column) and then crosses
  back across a $t$ edge in the same vertical strip.  Since the path
  in $\widetilde Y$ is immersed, the two $t$ edges must be distinct
  and the portion between them must contain $a^m$, $b^n$ or their
  inverses.  Actually this shows more than is claimed in the statement
  of the theorem.  Every reduced word equivalent to the identity in
  $G$ contains a syllable of the form $a^k$ where $k$ is a multiple of
  $m$ or $b^\ell$ where $\ell$ is a multiple of $n$.

  Finally, to prove (4) we use the projection of the closed curve to
  $T$ described above and systematically use the relation $a^m = b^n$
  to shrink the number of edges that the projection crosses in $T$.
  We should also note that with the lengths as assigned, this results
  in a nonlength increasing solution for the word problem of the torus
  knot group $G$.
\end{proof}

We conclude our discussion of torus knot groups by noting that in
addition to being easier to visualize, the geometry of $Y = Y_{p,q}$
is more closely tied to the geometry of the corresponding torus knot.

\begin{figure}
\includegraphics{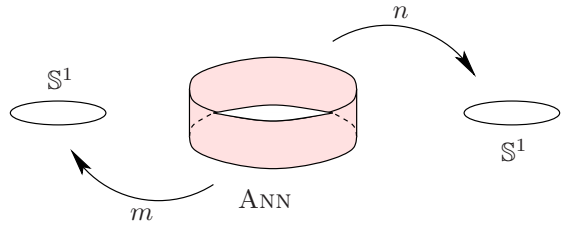}
\caption{If an annulus has its boundary cycles attached to two
  different circles with winding numbers $m$ and $n$, respectively,
  then the fundamental group of the resulting complex is the torus
  knot group $\tor(m,n)$.\label{fig:torus-annulus}}
\end{figure}

\begin{rem}[Torus knots and the complex $Y$]
  Let $p$ and $q$ be relatively prime integers.  There is a natural
  embedding of $Y = Y_{p,q}$ into the complement of the $(p,q)$-torus
  knot $K$ so that the complement of the knot deformation retracts
  onto $Y$.  We start by noting an alternate description of the space
  $Y$.  Imagine identifying the two edges of the polygon labeled $t$
  before performing the other edge identifications.  This shows that
  $Y$ can be also constructed by attaching an annulus to a pair of
  circles so that one boundary component is attached to one of the
  circles with winding number $m$ while attaching the other boundary
  component to the other circle with winding number $n$.  The result
  is homeomorphic to $Y$.  See Figure~\ref{fig:torus-annulus}.  To
  embed $Y$ into $\sph^3 \setminus K$ we sent the two circles to the
  core curves running through the centers of the two solid tori.  The
  annulus can then be embedded in $\sph^3$ and attached to the circles
  as needed to form $Y$ in such a way that it cuts through the
  boundary torus in a $(p,q)$-curve that is parallel to but disjoint
  from the original $(p,q)$-curve $K$.  Finally, it is not too
  difficult to construct an explicit deformation retraction from
  $\sph^3\setminus K$ to $Y$.
\end{rem}

%%%%%%%%%%%%%%%%%%%%%%%%%%%%%%%%%%%%%%%%%%%%%%%%%
\section{Solvable Baumslag-Solitar groups}
%%%%%%%%%%%%%%%%%%%%%%%%%%%%%%%%%%%%%%%%%%%%%%%%%

Our next family of examples are the solvable Baumslag-Solitar groups.
Althought these are groups where the standard $2$-complex adequately
encodes their geometry, we include a very brief discussion of their
basic properties so that we can can refer to them in the next section.

\begin{defn}[Baumslag-Solitar groups]
  The \emph{Baumslag-Solitar group} $BS(m,n)$ is the group defined by
  the presentation $\langle a,t \mid a^mt = ta^n \rangle$. The
  similarity between the Baumslag-Solitar group $BS(m,n)$ and the
  torus knot group $\tor(m,n)$ should be clear.  It too can be
  described as the fundamental group of a space obtained by attached
  the two boundary cycles to circles with winding numbers $m$ and $n$.
  The difference is that this time both boundary cycles are attached
  to the same circle. 
\end{defn}

A Baumslag-Solitar group is \emph{solvable} (in the classical sense of
that term) iff $m=n$ and this is also the case where the geometry is
most pleasing.  Let $G = BS(m,m)$ be a solvable Baumslag-Solitar group
and let $X$ be the standard $2$-complex of the presentation $\langle
a,t \mid a^m t = t a^m \rangle$.  The one polygon involved can be
given the structure of a rectangle as before and the universal cover
$\widetilde X$ has the structure of a tree $T$ cross the reals.  The
tree $T$ is a uniformly $m$-branching tree.  From this structure, one
can compute the center of $G$, see that $G$ is virtually
free-by-cyclic, solve the word problem in $G$ and establish basic
properties of reduced words equal to the identity in $G$.  In other
words, one can prove a theorem analogous to Theorem~\ref{thm:torus}.

%%%%%%%%%%%%%%%%%%%%%%%%%%%%%%%%%%%%%%
\section{One-relator Artin groups}
%%%%%%%%%%%%%%%%%%%%%%%%%%%%%%%%%%%%%%

Our third family of examples is closely related to the two previous
families.  Recall that an Artin group is defined by a presentation
inspired by Artin's classical presentation for the braid groups
\cite{Ar26,Ar47}.  In particular, they are defined by presentations in
which every relation is one of Artin's relations.

\begin{defn}[Artin relations and Artin groups]
  Let $(a,b)_m$ be the word of length $m$ which starts with $a$ and
  alternates between $a$ and $b$.  In symbols $(a,b)_m = abab\ldots$
  with $m$ letters total.  For example $(a,b)_2 = ab$, $(a,b)_3 = aba$
  and $(a,b)_4 = abab$.  An \emph{Artin relation} is a relation of the
  form $(a,b)_m = (b,a)_m$ with $m>1$.  Thus for small values of $m$
  we have commutation $ab=ba$, the braid relation $aba=bab$, and
  $abab=baba$ for $m=4$.  An \emph{Artin group} is a group defined by
  a presentation in which every relation is an Artin relation and
  there is at most one Artin relation for every distinct pair of
  generators.
\end{defn}

\begin{figure}
  \bt{cc}
  \bt{c}\includegraphics{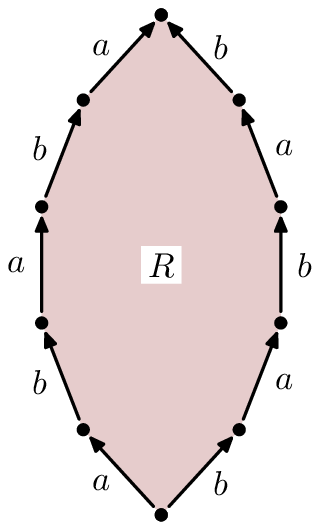}\et & 
  \bt{c}\includegraphics{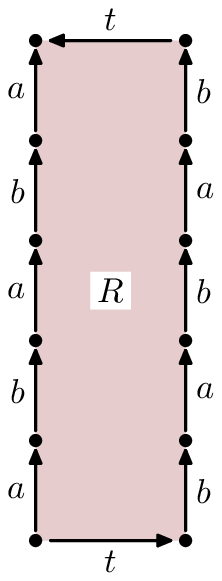}\et \\
  \bt{c}\includegraphics{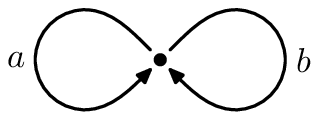}\et & 
  \bt{c}\includegraphics{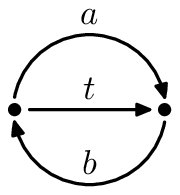}\et \\
  \et  
  \caption{The $2$-cell on the left attached to the graph on the left
    and the $2$-cell on the right attached to the graph on the right
    are both $2$-complexes with fundamental group $\art_5$.
    Contracting the edge $t$ shows that the two spaces are homotopy
    equivalent.\label{fig:artin-cells}}
\end{figure}

The simplest Artin groups are those with only two generators and one
relation.  Let $\art_m$ denote the \emph{one-relator Artin group}
defined by the presentation $\langle a,b \mid (ab)_m = (ba)_m \rangle$
and let $X = X_m$ be the corresponding standard $2$-complex.  As with
torus knots, the global structure of the universal cover of $X$ is not
immediately clear but there is a homotopy equivalent two-vertex
$2$-complex whose universal cover is much easier to visualize.  The
idea is to replace the Artin relation $(a,b)_m = (b,a)_m$ with a
similar relation that produces a $2$-complex whose $1$-skeleton looks
like the lower righthand side of Figure~\ref{fig:artin-cells}.  In
this graph it is possible to read the word $(a,b)_m$ and the word
$(b,a)_m$ without inserting any $t$ edges, but $t$ edges must be
inserted at the two transitions between the two words.  Also note that
the direction the $t$ edge needs to crossed depends on the parity of
$m$.  Explicitly, when $m$ is even we consider the combinatorial
description $[(a,b)_m t = t (b,a)_m ]$ and when $m$ is odd we consider
the combinatorial description $[(a,b)_m = t (b,a)_m t]$.  When these
relations are drawn as a rectangle similar to the one shown in
Figure~\ref{fig:artin-cells}, the two $t$ edges are pointing in the
same direction when $m$ is even and opposite directions when $m$ is
odd.

\begin{figure}
\includegraphics{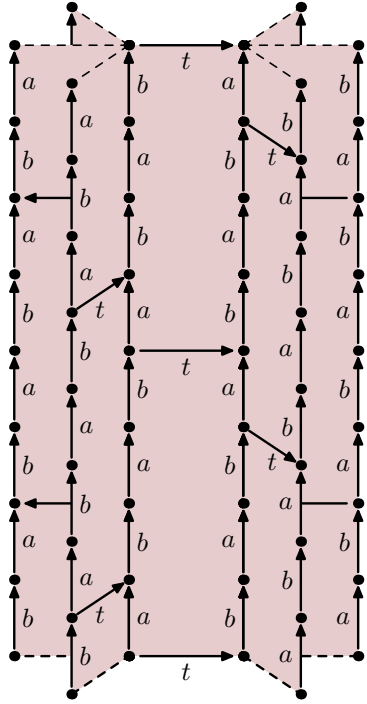}
\caption{A portion of the universal cover of $Y_4$.\label{fig:artin-u-cover}}
\end{figure}

In both cases, let $Y = Y_m$ denote the corresponding two-vertex
$2$-complex.  As in our previous examples, the one polygon involved
can be given a rectangular metric (with right angles at the endpoints
of the two $t$ edges) under which the universal cover $\widetilde Y$
is a metric direct product of a uniformly $m$-branching tree $T$ and
the reals.  A portion of the universal cover for for $Y_4$ is shown in
Figure~\ref{fig:artin-u-cover}.  From this structure, one can prove a
theorem analogous to Theorem~\ref{thm:torus}.

\begin{thm}[One-relator Artin groups]\label{thm:artin}
  If $G = \art_m$ is a one-relator Artin group with $m>1$ and $z$ is
  the element of $G$ represented by $(a,b)_m = (b,a)_m$ then: (1) the
  center of $G$ is an infinite cyclic subgroup generated by $z$ when
  $m$ is even and by $z^2$ when $m$ is odd; (2) $G$ is virtually a
  direct product of a free group and an infinite cyclic group; (3)
  every nontrivial reduced word equivalent to the identity in $G$
  contains a subword equal to $(a,b)_m$ or $(b,a)_m$ or their
  inverses; and finally, (4) every word equivalent to the identity in
  $G$ can be reduced to the identity by iteratively replacing
  $(a,b)_m$ with $(b,a)_m$, replacing $(b,a)_m$ with $(a,b)_m$, and
  performing free reductions.
\end{thm}

\begin{proof}[Proof Sketch]
  The proofs of these properties are nearly identical to the ones
  given for torus knots and they follow fairly quickly once the
  geometry of $\widetilde Y$ and action of $G$ is understood. First
  note that because contracting the edge labeled $t$ in $Y$ yields
  $X$, contracting the disjoint edges labeled $t$ in $\widetilde Y$
  yields $\widetilde X$.  Thus, we once again treat the edges labeled
  $t$ in $\widetilde Y$ as though they were contracted (without
  actually contracting them) allowing us to work with the
  geometrically pleasing $1$-skeleton of $\widetilde Y$ to establish
  results about the $1$-skeleton of $\widetilde X$, i.e. the Cayley
  graph of $G$ with respect to the generating set $\{a,b\}$.  For
  example, in $X$, the edges labeled $a$ and $b$ label loops which
  represent elements of the fundamental group, and as such they act on
  $\widetilde X$ by deck transformations once we have chosen a vertex
  in $\widetilde X$ as our base vertex.  Thus they also represent
  actions on $\widetilde Y$ by deck transformations once we have
  chosen an edge labeled $t$ as our base edge.  This time though, the
  deck transformations corresponding to $a$ and $b$ are more
  complicated.  The simple deck transformations are those associated
  with the words $ab$ and $ba$.  One of these stabilizes the column
  attached to the beginning of our $t$ edge, shifting it up two units
  while rotating the rest of $\widetilde Y$ around this column and the
  other performs a similar action with respect to the column attached
  to the other end.  After $m$ iterations of the first motion the
  entire complex experiences a pure vertical shift with no twisting.
  Similarly after $m$ iterations of the second motion.  Any element in
  the center of $G$ must commute with both of these motions and one
  can show that the possibilities are words that rigidly vertically
  shift the vertical strip containing our base $t$ edge.  When $m$ is
  odd, the smallest such shift is represented by $z^2 = (ab)^m =
  (ba)^m$ but when $m$ is even, the word $(a,b)_m$ representing $z$,
  equal to $(ab)^{m/2}$, also represents a rigid vertical shift and in
  both cases these elements are indeed central in $G$.

  To see (2) we construct is a finite-sheeted cover of $Y$ as before.
  First identify edges labeled $t$ in $\widetilde Y$ when they belong
  to the same vertical strip and are oriented in the same direction,
  identify two vertical strips if they have their $t$ edges at the
  same set of heights, and identify columns based on the parity of the
  heights at which the $a$ edges occur.  Geometrically the result is a
  direct product $\Gamma \times \sph^1$, where $\Gamma$ is a finite
  graph with two vertices and $m$ edges connecting them.  Thus its
  fundamental group $\F \times \Z$ where $\F$ is a free group of rank
  $m-1$ and since the cover is finite-sheeted, the subgroup this
  represents is finite index in $G$.

  For (3) we convert a reduced word equivalent to the identity into a
  closed immersed path in the $1$-skeleton of $\widetilde X$ starting
  at its base vertex and then finally to a closed immersed path in the
  $1$-skeleton of $\widetilde Y$ by traversing $t$ edges when
  necessary in order to continue.  Given such a path, we can look at
  its projection into the tree $T$.  The projection cannot be trivial
  since there are no closed immersed paths that remain in a single
  column and the projected curve cannot remain immersed since $T$ is a
  tree.  Thus there is a point in the projected curve where it crosses
  an edge of $T$ and then immediately backtracks across the same edge.
  If we consider the portion of the path in $\widetilde Y$ that
  produces this behavior, we see a path that crosses a $t$ edge,
  travels up or down a column and then crosses back across a $t$ edge
  in the same vertical strip.  Since the path in $\widetilde Y$ is
  immersed, the two $t$ edges must be distinct and the portion between
  them must contain $(a,b)_m$, $(b,a)_m$ or their inverses.

  Finally, to prove (4) we use the projection of the closed curve to
  $T$ described above and systematically use the relation $(a,b)_m =
  (b,a)_m$ to shrink the number of edges that the projection crosses
  in $T$ and note that this results in a nonlength increasing solution
  for the word problem of $G$.
\end{proof}

The fact that one-relator Artin groups have properties similar to
torus knot groups and solvable Baumslag-Solitar groups is not
accidental.

\begin{figure}
\includegraphics{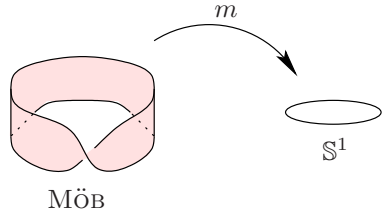}
\caption{A M\"obius band attached to a circle.\label{fig:moebius}}
\end{figure}

\begin{rem}[Relations with previous examples]
  Let $G$ be the one-relator Artin group $\art_m$ and let $Y$ be the
  corresponding two-vertex complex described above with fundamental
  group $G$.  If we identify the two edges labeled $t$ before carrying
  out the other identifications, we can reanalyze $Y$, depending on
  the parity of $m$, as either a solvable Baumslag-Solitar group or a
  torus knot group.  More specifically, when $m$ is even, identifying
  the two $t$ edges creates a torus with boundary cycles attached to
  the same circle.  In particular, the Artin group $\art_m$ is
  isomorphic to solvable Baumslag-Solitar group $BS(m/2,m/2)$ when $m$
  is even.  On the other hand, when $m$ is odd, identifying the two
  $t$ edges creates a m\"obius strip whose boundary cycle is attached
  to a circle with winding number $m$ (Figure~\ref{fig:moebius}).
  Cutting the m\"obius strip along its central curve shows that
  $\art_m$ is isomorphic to the torus knot group $\tor(2,m)$ when $m$
  is odd.
\end{rem}

Note that from a topological perspective, all three classes of groups
(torus knot groups, Baumslag-Solitar groups, one-relator Artin groups)
are extremely simple as they only involve attaching annuli or m\"obius
strips to one or more circles.  It is thus somewhat surprising that
they are treated separately and that there does not exist a more
uniform set of notations for these groups and their $2$-complexes.
Finally, we conclude as promised with a short geometric proof of a key
lemma used by Kenneth Appel and Paul Schupp in their investigation of
Artin groups of large and extra-large type.

\begin{defn}[Large type and extra large type]
  As mentioned earlier, an \emph{Artin group} is a group $G$ defined
  by a presentation in which every relation is an Artin relation
  $(a,b)_m = (b,a)_m$ and for every pair of generators there is at
  most one such relation.  If for every relation in the presentation,
  the integer $m$ is at least $3$ then $G$ is an Artin group of
  \emph{large type} and if every integer $m$ is at least $4$ then $G$
  is an Artin group of \emph{extra large type}.
\end{defn}

In order to analyze large and extra-large Artin groups using small
cancellation theory, they first needed to prove a key intermediate
result specifically about one-relator Artin groups.  In particular,
Appel and Schupp used a detailed inductive analysis of possible van
Kampen diagrams in order to establish the follow result that occurs as
Lemma~$6$ in \cite{ApSch83}.

\begin{prop}[Syllable counts]\label{prop:syllables}
  Every nontrivial cyclically reduced word that is equal to the
  identity in $\art_m$ contains at least $2m$ syllables.
\end{prop}

\begin{proof}
  Let $G$ be $\art_m$, let $X$ be standard $2$-complex and let $Y$ be
  the two-vertex complex with fundamental group $G$ described above.
  The proof proceeds by counting breaks between syllables using the
  vertical projection from $Y$ to the tree $T$ and the horizontal
  projection from $Y$ to the real line $\R$.  As in the proof of
  Theorem~\ref{thm:artin}, we lift the nontrivial cyclically reduced
  word equal to the identity in $G$ to an immersed loop in the
  universal cover $\widetilde X$ and then to an immersed loop in the
  universal cover $\widetilde Y$.  As argued above, the projection to
  $T$ is nontrivial and no longer immersed.  Moreover, each time the
  projection to $T$ crosses an edge and immediately recrosses it in
  the other direction, we can find a copy of the word $(a,b)_m$ or
  $(b,a)_m$ in $\widetilde Y$ and each such occurence includes $(m-1)$
  syllable breaks, i.e. gaps between letters with distinct generators
  on either side.  Since any nontrivial path in a tree includes at
  least two such backtracks, we have found $2(m-1) = 2m-2$ syllable
  breaks.  The final two syllable breaks are located by projecting
  horizontally to the real line.  The path in $\widetilde Y$ projects
  to a path moving up and down the real line.  Because it is a closed
  loop, it attains a maximum and a minimum value.  When the path
  reaches its maximum value it must change columns in $\widetilde Y$
  (i.e. it must cross a $t$ edge) before continuing back down since
  the path in $\widetilde Y$ is immersed.  This leads to a subword of
  the form $aB$ or $bA$ and to a syllable break that was not
  previously counted.  Similarly local minima in the projection lead
  to subwords of the form $Ab$ or $Ba$ and to a final syllable break
  that was not previous counted.  Because the closed path has at least
  $2(m-1) + 2 = 2m$ syllable breaks, it must contain at least $2m$
  syllables.
\end{proof}

In \cite{ApSch83} Appel and Schupp used this result to analyze Artin
groups of extra-large type and in \cite{Ap84} Appel extended the
analysis to Artin groups of large type.  Roughly speaking, if you
consider van Kampen diagrams over Artin groups of with respect to an
infinite presentation that includes every nontrivial cyclically
reduced word in a subgroup generated by two of its generators, then
the one can alway find a diagram where no two cells sharing an edge
have boundary cycles from the same two generator subgroup.  Under
these conditions, the overlap between two cells consists of a single
letter and thus lives within a single syllable of the boundary word of
the either cell.  The large or extra-large condition, combined with
Proposition~\ref{prop:syllables} means that these diagrams satisfy the
small cancellation conditions $C(6)$ or $C(8)$, respectively.  Once
the tools of small cancellation theory are available, they are then
able to establish many foundational results for large and extra-large
Artin groups.  In their original paper Appel and Schupp needed to work
a bit to establish Proposition~\ref{prop:syllables}.  The geometry of
$Y$ makes this proposition much more transparent.

\bibliographystyle{plain}

%%%%%%%%%%%%%%
\end{document}